\documentclass[11pt,reqno]{amsart}
\usepackage{amsmath,amssymb}
\usepackage[colorlinks=true,linkcolor=magenta,citecolor=blue]{hyperref}
\usepackage[margin=1in]{geometry}
\usepackage{enumitem} 
\usepackage{mathabx}
\usepackage{epsfig}

\newtheorem{dn}{Definition}[section]
\newtheorem{dl}{Theorem}[section]
\newtheorem{md}{Proposition}[section]
\newtheorem{bd}{Lemma}[section]
\newtheorem{hq}{Corollary}[section]
\newtheorem{nx}{Remark}[section]
\newtheorem{vd}{Example}[section]
\newcommand{\R}{\mathbb{R}}

\newcommand{\Z}{\mathbb{Z}}

\newcommand{\e}{\varepsilon}
\newcommand{\ity}{\infty}
\newcommand{\f}{\frac}

\newcommand{\bbd}{\begin{bd}}
\newcommand{\ebd}{\end{bd}}
\newcommand{\bdn}{\begin{dn}}
\newcommand{\edn}{\end{dn}}
\newcommand{\bhq}{\begin{hq}}
\newcommand{\ehq}{\end{hq}}
\newcommand{\bdl}{\begin{dl}}
\newcommand{\edl}{\end{dl}}
\newcommand{\bnx}{\begin{nx}}
\newcommand{\enx}{\end{nx}}
\newcommand{\bmd}{\begin{md}}
\newcommand{\emd}{\end{md}}
\newcommand{\bvd}{\begin{vd}}
\newcommand{\evd}{\end{vd}}

\title[Critical exponents for a structurally damped wave system]{Existence and nonexistence of global solutions for a structurally damped wave system with power nonlinearities}
\author{Tuan Anh Dao}
\address{Tuan Anh Dao \hfill\break
$\quad$ School of Applied Mathematics and Informatics, Hanoi University of Science and Technology, No.1 Dai Co Viet road, Hanoi, Vietnam \hfill\break
Faculty for Mathematics and Computer Science, TU Bergakademie Freiberg, Pr\"{u}ferstr. 9, 09596, Freiberg, Germany}
\email{anh.daotuan@hust.edu.vn}

\begin{document}

	
\begin{abstract}
Our interest itself of this paper is strongly inspired from an open problem in the paper \cite{Dabbicco} published by D'Abbicco. In this article, we would like to study the Cauchy problem for a weakly coupled system of semi-linear structurally damped wave equations. Main goal is to find the threshold, which classifies the global (in time) existence of small data solutions or the nonexistence of global solutions under the growth condition of the nonlinearities.
\end{abstract}

\maketitle

\section{Introduction} \label{Sec.main}
In this paper, let us consider the following Cauchy problem for weakly coupled system of semi-linear structurally damped wave equations:
\begin{equation}
\begin{cases}
u_{tt}- \Delta u+ (-\Delta)^{\delta_1} u_t=|v|^{p}, &\quad x\in \R^n,\, t \ge 0, \\
v_{tt}- \Delta v+ (-\Delta)^{\delta_2} v_t=|u|^{q}, &\quad x\in \R^n,\, t \ge 0, \\
u(0,x)= u_0(x),\quad u_t(0,x)=u_1(x), &\quad x\in \R^n, \\
v(0,x)= v_0(x),\quad v_t(0,x)=v_1(x), &\quad x\in \R^n, \label{pt1.1}
\end{cases}
\end{equation}
for any $\delta_1,\,\delta_2 \in [0,1]$ and for nonlinearities with powers $p,\, q >1$. The special case of (\ref{pt1.1}) with $\delta_1= \delta_2= \frac{1}{2}$ in the form
\begin{equation} \label{pt1.2}
\begin{cases}
u_{tt}- \Delta u+ (-\Delta)^{\frac{1}{2}} u_t= |v|^p,\quad  v_{tt}- \Delta v+ (-\Delta)^{\frac{1}{2}} v_t= |u|^q, \\
u(0,x)= u_0(x),\quad u_t(0,x)=u_1(x),\quad v(0,x)= v_0(x),\quad v_t(0,x)=v_1(x),
\end{cases}
\end{equation}
was well- studied by D'Abbicco in \cite{Dabbicco}. In the cited paper, he succeeded to determine the critical exponent for (\ref{pt1.2}). For details, the author proved the global (in time) existence of small data solutions to (\ref{pt1.2}) in any space dimensions $n\ge 2$ if the condition
$$ \frac{1+ \max\{p,\,q\}}{pq-1}< \frac{n-1}{2} $$
holds by using sharp decay estimates for solutions to the linear corresponding Cauchy problem. Moreover, the above condition is sharp because a nonexistence result of global (in time) weak solutions to (\ref{pt1.2}) was also discussed if this condition is no longer true. The proof of blow-up result is based on a contradiction argument by using the test function method (see, for example, \cite{Dabbicco,Zhang}). The fact is that for this purpose some difficulties arise. In general, standard test function method, i.e. test functions with compact support, is not directly applicable since this method relies on pointwise control of derivatives of test functions. In addition, the fractional Laplacian operators $(-\Delta)^\delta$ for any $\delta \in (0,1)$ are well-known non-local operators, it follows that supp$(-\Delta)^\delta \phi$ is bigger than supp$\phi$ for any $\phi \in \mathcal{C}_0^\ity(\R^n)$ in general. However, this application linked to the estimate
$$ (-\Delta)^\delta \phi^\ell \le \ell \phi^{\ell-1}(-\Delta)^\delta \phi \quad \text{ for }\delta \in (0,1),\, \ell \ge 1 \text{ and for all } \phi\ge 0, \phi \in \mathcal{C}_0^\ity(\R^n) $$
is possible to (\ref{pt1.2}) due to the following key observation: Any local or global solution to (\ref{pt1.2}) is nonnegative with the assumption of nonnegative initial data $u_1,v_1$ and $u_0=v_0=0$, which was investigated by D'Abbicco-Reissig in \cite{DabbiccoReissig}. Thanks to this essential property, the above inequality works well to extend the test function method to (\ref{pt1.2}). Unfortunately, we cannot expect nonnegative solutions to (\ref{pt1.1}), which contains the nonlocal terms $(-\Delta)^{\delta}$ for any $\delta \in (0,1)$.
\par For this reason, the first main motivation of this paper is to prove the global (in time) existence of small data solutions to (\ref{pt1.1}), where the parameters $\delta_1$ and $\delta_2$ are not necessary to be equal. More in details, we would like to explain the impact of the flexible choice of the parameters $\delta_1,\,\delta_2$ on our global (in time) existence results and the range of admissible exponents $p,\,q$ as well. To establish this, we have in mind to take advantage of the better decay estimates available for the corresponding linear wave equations with structural damping $(-\Delta)^{\delta}u_t$ of (\ref{pt1.1}) in the following form:
\begin{equation}
w_{tt}- \Delta w+ (-\Delta)^{\delta} w_t=0,\quad w(0,x)= w_0(x),\quad w_t(0,x)= w_1(x), \label{pt1.3}
\end{equation}
where $\delta= \delta_1$ or $\delta= \delta_2$. From these appearing difficulties as mentioned above, the second main motivation of this paper is to find the precise critical exponents to (\ref{pt1.1}) with general cases of $\delta_1,\,\delta_2 \in [0,1]$, especially we are interested in facing up to the proof of blow-up result, where the requirement of nonnegativity of solutions does not appear for (\ref{pt1.1}). In order to overcome this difficulty, the crux of our ideas is to apply a modified test function method effectively in dealing with the fractional Laplacian $(-\Delta)^{\delta_1}$ and $(-\Delta)^{\delta_2}$.
\par Moreover, concerning the linear equation (\ref{pt1.3}) and some of its semi-linear equations with the power nonlinearity $|u|^p$ we want to point out the paper \cite{DabbiccoReissig} of D'Abbicco-Reissig. The authors have proposed to distinguish between ``parabolic like models" in the case $\delta \in [0,\frac{1}{2}]$, the so-called effective damping, and ``hyperbolic like models" or ``wave like models" in the case $\delta \in (\frac{1}{2},1]$, the so-called noneffective damping. To the best of author's knowledge, it seems that nobody has ever succeeded to determine really critical exponent to semi-linear structurally damped wave equations with noneffective damping. Hence, it is still an open problem as far as to explore. From this observation, in order to give a partial positive answer to the open problem in \cite{Dabbicco}, it is quite natural that we may restrict ourselves to consider only (\ref{pt1.1}) with effective damping, i.e. the assumption of $\delta _1,\, \delta _2 \in [0,\frac{1}{2}]$ is of our interest in this paper.

\subsection{Notations}
We use the following notations throughout this paper.
\begin{itemize}[leftmargin=*]
\item We write $f\lesssim g$ when there exists a constant $C>0$ such that $f\le Cg$, and $f \approx g$ when $g\lesssim f\lesssim g$.
\item As usual, the spaces $H^a$ and $\dot{H}^a$ with $a \ge 0$ stand for Bessel and Riesz potential spaces based on $L^2$ spaces. Here $\big<D\big>^a$ and $|D|^a$ denote the pseudo-differential operators with symbols $\big<\xi\big>^a$ and $|\xi|^a$, respectively. We denote $\widehat{f}(t,\xi):= \mathfrak{F}_{x\rightarrow \xi}\big(f(t,x)\big)$ as the Fourier transform with respect to the space variable of a function $f(t,x)$. 
\item For a given number $s \in \R$, we denote
$$ [s]:= \max \big\{k \in \Z \,\, : \,\, k\le s \big\} \quad \text{ and }\quad [s]^+:= \max\{s,0\} $$
as its integer part and its positive part, respectively.
\item We put $\big< x\big>:= \sqrt{1+|x|^2}$, the so-called Japanese bracket of $x \in \R^n$.
\item We fix the constant $m_0:=\frac{2m}{2-m}$, that is, $\frac{1}{m_0}=\frac{1}{m}- \frac{1}{2}$ with $m \in [1,2)$.
\item Finally, we introduce the spaces
$\mathcal{A}:= \big(L^m \cap H^1\big) \times \big(L^m \cap L^2\big)$ with the norm
$$\|(u_0,u_1)\|_{\mathcal{A}}:=\|u_0\|_{L^m}+ \|u_0\|_{H^1}+ \|u_1\|_{L^m}+ \|u_1\|_{L^2}, \quad \text{ where }m \in [1,2). $$
\end{itemize}

\subsection{Main results}
Let us state the main results which will be proved in the present paper.
\bdl[\textbf{Global existence for $\delta_1 \ge \delta_2$}] \label{dl1.1}
Let us assume $\delta _1,\, \delta _2 \in \big[0,\frac{1}{2}\big]$ and $\delta_1 \ge \delta_2$. Let $m\in [1,2)$ and $n> 2m_0\delta_1$. We assume that the conditions are satisfied
\begin{align}
&\frac{2}{m} \le p,\, q < \ity &\text{ if }&\quad n \le 2, \label{GN11A1} \\
&\frac{2}{m} \le p,\, q \le \frac{n}{n- 2} &\text{ if }&\quad 2 < n \le \frac{4}{2-m}. \label{GN11A2}
\end{align}
Moreover, we suppose the following conditions:
\begin{equation} \label{exponent11A1}
\frac{1+ q\frac{1-\delta_2}{1-\delta_1}+ (pq-1)\delta_2}{(q-1)\frac{\delta_1- \delta_2}{1- \delta_2}+ pq-1}< \frac{n}{2m},
\end{equation}
and
\begin{equation} \label{exponent11A2}
p \le 1+ \frac{2m}{n- 2m\delta_2} \le 1+ \frac{2m}{n- 2m\delta_1} < q.
\end{equation}
Then, there exists a constant $\e_0>0$ such that for any small data
$$ \big((u_0,u_1),\, (v_0,v_1) \big) \in \mathcal{A} \times \mathcal{A} \text{ satisfying the assumption } \|(u_0,u_1)\|_{\mathcal{A}}+ \|(v_0,v_1)\|_{\mathcal{A}} \le \e_0, $$
we have a uniquely determined global (in time) small data energy solution
$$ (u,v) \in \Big(C\big([0,\ity),H^1\big)\cap C^1\big([0,\ity),L^2\big)\Big)^2 $$
to (\ref{pt1.1}). The following estimates hold:
\begin{align}
\|u(t,\cdot)\|_{L^2} &\lesssim (1+t)^{-\frac{n}{2(1-\delta_1)}(\frac{1}{m}-\frac{1}{2})+ \frac{\delta_1}{1-\delta_1}+ \e(p,\delta_2)} \big(\|(u_0,u_1)\|_{\mathcal{A}}+ \|(v_0,v_1)\|_{\mathcal{A}}\big), \label{decayrate11A1} \\
\big\|\nabla u(t,\cdot)\big\|_{L^2} &\lesssim (1+t)^{-\frac{n}{2(1-\delta_1)}(\frac{1}{m}-\frac{1}{2})- \frac{1-2\delta_1}{2(1-\delta_1)}+ \e(p,\delta_2)} \big(\|(u_0,u_1)\|_{\mathcal{A}}+ \|(v_0,v_1)\|_{\mathcal{A}}\big), \label{decayrate11A2} \\
\|u_t(t,\cdot)\|_{L^2} &\lesssim (1+t)^{-\frac{n}{2(1-\delta_1)}(\frac{1}{m}-\frac{1}{2})- \frac{1-2\delta_1}{1-\delta_1}+ \e(p,\delta_2)} \big(\|(u_0,u_1)\|_{\mathcal{A}}+ \|(v_0,v_1)\|_{\mathcal{A}}\big), \label{decayrate11A3} \\
\|v(t,\cdot)\|_{L^2} &\lesssim (1+t)^{-\frac{n}{2(1-\delta_2)}(\frac{1}{m}-\frac{1}{2})+ \frac{\delta_2}{1-\delta_2}} \big(\|(u_0,u_1)\|_{\mathcal{A}}+ \|(v_0,v_1)\|_{\mathcal{A}}\big), \label{decayrate11A4} \\
\big\|\nabla v(t,\cdot)\big\|_{L^2} &\lesssim (1+t)^{-\frac{n}{2(1-\delta_2)}(\frac{1}{m}-\frac{1}{2})- \frac{1-2\delta_2}{2(1-\delta_2)}} \big(\|(u_0,u_1)\|_{\mathcal{A}}+ \|(v_0,v_1)\|_{\mathcal{A}}\big), \label{decayrate11A5} \\
\|v_t(t,\cdot)\|_{L^2} &\lesssim (1+t)^{-\frac{n}{2(1-\delta_2)}(\frac{1}{m}-\frac{1}{2})- \frac{1-2\delta_2}{1-\delta_2}} \big(\|(u_0,u_1)\|_{\mathcal{A}}+ \|(v_0,v_1)\|_{\mathcal{A}}\big), \label{decayrate11A6}
\end{align}
where $\e(p,\delta_2):= 1- \frac{n}{2m(1-\delta_2)}(p-1)+ \frac{p\delta_2}{1-\delta_2}+\e$ with a sufficiently small positive number $\e$.
\edl

\bdl[\textbf{Global existence for $\delta_2 \ge \delta_1$}] \label{dl1.2}
Let us assume $\delta _1,\, \delta _2 \in \big[0,\frac{1}{2}\big]$ and $\delta_2 \ge \delta_1$. Let $m\in [1,2)$ and $n> 2m_0\delta_2$. We assume that the conditions (\ref{GN11A1}) and (\ref{GN11A2}) are satisfied. Moreover, we suppose the following conditions:
\begin{equation} \label{exponent12A1}
\frac{1+ p\frac{1-\delta_1}{1-\delta_2}+ (pq-1)\delta_1}{(p-1)\frac{\delta_2- \delta_1}{1- \delta_1}+ pq-1}< \frac{n}{2m},
\end{equation}
and
\begin{equation} \label{exponent12A2}
q \le 1+ \frac{2m}{n- 2m\delta_1} \le 1+ \frac{2m}{n- 2m\delta_2} < p.
\end{equation}
Then, we have the same conclusions as in Theorem \ref{dl1.1}. But the estimates (\ref{decayrate11A1})-(\ref{decayrate11A6}) are modified in the following way:
\begin{align}
\|u(t,\cdot)\|_{L^2} &\lesssim (1+t)^{-\frac{n}{2(1-\delta_1)}(\frac{1}{m}-\frac{1}{2})+ \frac{\delta_1}{1-\delta_1}} \big(\|(u_0,u_1)\|_{\mathcal{A}}+ \|(v_0,v_1)\|_{\mathcal{A}}\big), \label{decayrate12A1} \\
\big\|\nabla u(t,\cdot)\big\|_{L^2} &\lesssim (1+t)^{-\frac{n}{2(1-\delta_1)}(\frac{1}{m}-\frac{1}{2})- \frac{1-2\delta_1}{2(1-\delta_1)}} \big(\|(u_0,u_1)\|_{\mathcal{A}}+ \|(v_0,v_1)\|_{\mathcal{A}}\big), \label{decayrate12A2} \\
\|u_t(t,\cdot)\|_{L^2} &\lesssim (1+t)^{-\frac{n}{2(1-\delta_1)}(\frac{1}{m}-\frac{1}{2})- \frac{1-2\delta_1}{1-\delta_1}} \big(\|(u_0,u_1)\|_{\mathcal{A}}+ \|(v_0,v_1)\|_{\mathcal{A}}\big), \label{decayrate12A3} \\
\|v(t,\cdot)\|_{L^2} &\lesssim (1+t)^{-\frac{n}{2(1-\delta_2)}(\frac{1}{m}-\frac{1}{2})+ \frac{\delta_2}{1-\delta_2}+ \e(q,\delta_1)} \big(\|(u_0,u_1)\|_{\mathcal{A}}+ \|(v_0,v_1)\|_{\mathcal{A}}\big), \label{decayrate12A4} \\
\big\|\nabla v(t,\cdot)\big\|_{L^2} &\lesssim (1+t)^{-\frac{n}{2(1-\delta_2)}(\frac{1}{m}-\frac{1}{2})- \frac{1-2\delta_2}{2(1-\delta_2)}+ \e(q,\delta_1)} \big(\|(u_0,u_1)\|_{\mathcal{A}}+ \|(v_0,v_1)\|_{\mathcal{A}}\big), \label{decayrate12A5} \\
\|v_t(t,\cdot)\|_{L^2} &\lesssim (1+t)^{-\frac{n}{2(1-\delta_2)}(\frac{1}{m}-\frac{1}{2})- \frac{1-2\delta_2}{1-\delta_2}+ \e(q,\delta_1)} \big(\|(u_0,u_1)\|_{\mathcal{A}}+ \|(v_0,v_1)\|_{\mathcal{A}}\big), \label{decayrate12A6}
\end{align}
where $\e(q,\delta_1):= 1- \frac{n}{2m(1-\delta_1)}(q-1)+ \frac{q\delta_1}{1-\delta_1}+\e$ with a sufficiently small positive number $\e$.
\edl

\begin{nx}
\fontshape{n}
\selectfont
Here we want to stress out that $\e(p,\delta_2)$ and $\e(q,\delta_1)$ appearing in Theorems \ref{dl1.1} and \ref{dl1.2} represent some loss of decay in comparison with the corresponding decay estimates for solutions to (\ref{pt1.3}) (see later, Corollary \ref{hq2.1}). Besides, thanks to the conditions (\ref{exponent11A1}) and (\ref{exponent12A1}), both $\e(p,\delta_2)$ and $\e(q,\delta_1)$ are nonnegative.
\end{nx}

Finally, in order to show the optimality of our exponents to (\ref{pt1.1}), we have the following blow-up results.

\bdl[\textbf{Blow-up for initial data in $L^1$}] \label{dl1.3}
Let $\delta_1,\,\delta_2 \in \big[0,\frac{1}{2}\big]$. We assume that we choose the initial data $u_0=v_0=0$ and $u_1,\,v_1 \in L^1$ satisfying the following relations:
\begin{equation} \label{optimal13.1}
\int_{\R^n} u_1(x)dx > \epsilon_1 \quad \text{ and }\quad \int_{\R^n} v_1(x)dx > \epsilon_2,
\end{equation}
where $\epsilon_1$ and $\epsilon_2$ are suitable nonnegative constants. Moreover, we suppose the following conditions:
\begin{equation}
\frac{n}{2}\le \frac{1+ q\frac{1-\delta_2}{1-\delta_1}+ (pq-1)\delta_2}{(q-1)\frac{\delta_1- \delta_2}{1- \delta_2}+ pq-1} \quad \text{ if }\quad \delta_1 \ge \delta_2, \label{optimal13.2}
\end{equation}
or
\begin{equation}
\frac{n}{2}\le \frac{1+ p\frac{1-\delta_1}{1-\delta_2}+ (pq-1)\delta_1}{(p-1)\frac{\delta_2- \delta_1}{1- \delta_1}+ pq-1} \quad \text{ if }\quad \delta_2 \ge \delta_1. \label{optimal13.3}
\end{equation}
Then, there is no global (in time) Sobolev solution $(u,v) \in C\big([0,\infty),L^2\big) \times C\big([0,\infty),L^2\big)$ to (\ref{pt1.1}).
\edl

\begin{nx}
\fontshape{n}
\selectfont
If we choose $m=1$ in Theorems \ref{dl1.1} and \ref{dl1.2}, then from Theorem \ref{dl1.3} it follows that the exponents $p,\,q$ given by
$$ \frac{1+ q\frac{1-\delta_2}{1-\delta_1}+ (pq-1)\delta_2}{(q-1)\frac{\delta_1- \delta_2}{1- \delta_2}+ pq-1}= \frac{n}{2} \quad \text{ if }\quad \delta_1 \ge \delta_2, $$
or
$$ \frac{1+ p\frac{1-\delta_1}{1-\delta_2}+ (pq-1)\delta_1}{(p-1)\frac{\delta_2- \delta_1}{1- \delta_1}+ pq-1}= \frac{n}{2} \quad \text{ if }\quad \delta_2 \ge \delta_1 $$
are precisely critical in the case $\delta_1,\,\delta_2 \in \big[0,\frac{1}{2}\big]$.
\end{nx}
\begin{nx}
\fontshape{n}
\selectfont
We can observe that by setting $\delta_1=\delta_2=0$ or $\delta_1=\delta_2=\frac{1}{2}$, our main results from Theorems \ref{dl1.1}, \ref{dl1.2}, \ref{dl1.3} coincide really with those derived from \cite{NishiharaWakasugi} or \cite{Dabbicco}, respectively.
\end{nx}

\bdl[\textbf{Blow-up for initial data not in $L^1$}] \label{dl1.4}
Let $\delta_1= \delta_2=: \delta \in \big[0,\frac{1}{2}\big]$ and $m \in (1,2)$. We assume that we choose the initial data $u_0=v_0=0$ and $u_1,\,v_1 \in L^m$ satisfying the following relations:
\begin{equation} \label{optimal14.1}
u_1(x) \ge \epsilon_0 (1+|x|)^{-\frac{n+\epsilon_1}{m}} \quad \text{ and }\quad v_1(x) \ge \epsilon_0 (1+|x|)^{-\frac{n+\epsilon_2}{m}},
\end{equation}
where $\epsilon_0,\,\epsilon_1,\,\epsilon_2$ are suitably small positive constants. Moreover, we suppose the condition
\begin{equation}
\frac{n-2m\delta}{2m}< \frac{1+ \max\{p,\,q\}}{pq-1}. \label{optimal14.2}
\end{equation}
Then, there is no global (in time) Sobolev solution $(u,v) \in C\big([0,\infty),L^2\big) \times C\big([0,\infty),L^2\big)$ to (\ref{pt1.1}).
\end{dl}

\begin{nx}
\fontshape{n}
\selectfont
By plugging $\delta_1= \delta_2=: \delta$ into Theorems \ref{dl1.1} and \ref{dl1.2}, it is clear that we have found the critical exponents $p,\,q$ from Theorem \ref{dl1.4}. It only remains an open problem to verify whether there exist global (in time) Sobolev solutions or not in the following critical values:
$$ \frac{1+ \max\{p,\,q\}}{pq-1}= \frac{n-2m\delta}{2m}. $$
\end{nx}

\textbf{The outline of this article is presented as follows}: Section \ref{Sec.Pre} is to provide $(L^m \cap L^2)- L^2$ estimates and $L^2- L^2$ estimates for solutions to (\ref{pt1.2}), with $m \in [1,2)$, and some of essential properties of a modified test function method from the recent papers of Dao \cite{Dao1} and Dao-Reissig \cite{DaoReissig}, respectively. In Section \ref{Sec.Global existence}, we prove the global (in time) existence of small data solutions to (\ref{pt1.1}). Finally, we devote to the proof of nonexistence result of global solutions to (\ref{pt1.1}) in Section \ref{Sec.Blow-up}.

\section{Preliminaries} \label{Sec.Pre}
In this section, we collect some preliminary knowledge needed in our proofs.

\subsection{Linear estimates}
Main purpose is to recall $(L^m \cap L^2)- L^2$ and $L^2- L^2$ estimates for solutions and some of their derivatives to (\ref{pt1.3}) from the recent paper of Dao \cite{Dao1}. Using partial Fourier transformation to (\ref{pt1.3}) we have the following Cauchy problem:
\begin{equation}
\widehat{w}_{tt}+ |\xi|^{2\delta} \widehat{w}_t+ |\xi|^2 \widehat{w}=0,\quad \widehat{w}(0,\xi)= \widehat{w}_0(\xi),\quad \widehat{w}_t(0,\xi)= \widehat{w}_1(\xi). \label{pt2.1}
\end{equation}
The characteristic roots are
$$ \lambda_{1,2}=\lambda_{1,2}(\xi)= \f{1}{2}\Big(-|\xi|^{2\delta}\pm \sqrt{|\xi|^{4\delta}-4|\xi|^2}\Big). $$
The solutions to (\ref{pt2.1}) are written by the following form (here we assume $\lambda_{1}\neq \lambda_{2}$):
\begin{align*}
\widehat{w}(t,\xi)&= \frac{\lambda_1 e^{\lambda_2 t}-\lambda_2 e^{\lambda_1 t}}{\lambda_1- \lambda_2}\widehat{w}_0(\xi)+ \frac{e^{\lambda_1 t}-e^{\lambda_2 t}}{\lambda_1- \lambda_2}\widehat{w}_1(\xi) \\
&=: \widehat{K}_{0,\delta}(t,\xi)\widehat{w}_0(\xi)+ \widehat{K}_{1,\delta}(t,\xi)\widehat{w}_1(\xi).
\end{align*}
For this reason, we may read the solutions to (\ref{pt1.3}) as follows:
$$ w(t,x)=K_{0,\delta}(t,x) \ast_{x} w_0(x)+ K_{1,\delta}(t,x) \ast_{x} w_1(x). $$

\bmd[Proposition 2.3 in \cite{Dao1} with $\sigma=1$] \label{md2.1}
Let $\delta \in \big[0,\f{1}{2}\big]$ in (\ref{pt1.3}) and $m \in [1,2)$. Let $k\ge 0$, $j=0,1$. The solutions to (\ref{pt1.3}) satisfy the $(L^m \cap L^2)-L^2$ estimates
\begin{align*}
\big\|\partial_t^j \nabla^k w(t,\cdot)\big\|_{L^2} &\lesssim (1+t)^{-\frac{n}{2(1-\delta)}(\frac{1}{m}-\frac{1}{2})- \frac{k+2j\delta}{2(1-\delta)}}\|w_0\|_{L^m \cap H^{k+j}} \\ 
&\qquad \qquad + (1+t)^{-\frac{n}{2(1-\delta)}(\frac{1}{m}-\frac{1}{2})- \frac{k}{2(1-\delta)}-j+1}\|w_1\|_{L^m \cap H^{[k+(j-1)]^+}}
\end{align*}
and the $L^2-L^2$ estimates
$$ \big\|\partial_t^j \nabla^k w(t,\cdot)\big\|_{L^2} \lesssim (1+t)^{-\frac{k+2j\delta}{2(1-\delta)}}\|w_0\|_{H^{k+j}}+ (1+t)^{-\frac{k}{2(1-\delta)}-j+1}\|w_1\|_{H^{[k+(j-1)]^+}} $$
for all space dimensions $n \ge 1$.
\emd
\begin{nx}
\fontshape{n}
\selectfont
Here we want to underline that although all the decay estimates from Proposition \ref{md2.1} are available for any space dimensions $n\ge 1$, under a constraint condition to space dimensions $n> 2m_0\delta$ we may conclude the better decay estimates. Namely, we obtain the following result.
\end{nx}
\bmd[Proposition 2.2 in \cite{Dao1} with $\sigma=1$] \label{md2.2}
Let $\delta \in \big[0,\f{1}{2}\big]$ in (\ref{pt1.3}) and $m \in [1,2)$. Let $k\ge 0$, $j=0,1$. The solutions to (\ref{pt1.3}) satisfy the $(L^m \cap L^2)-L^2$ estimates
\begin{align*}
\big\|\partial_t^j \nabla^k w(t,\cdot)\big\|_{L^2} &\lesssim (1+t)^{-\frac{n}{2(1-\delta)}(\frac{1}{m}-\frac{1}{2})- \frac{k}{2(1-\delta)}-j}\|w_0\|_{L^m \cap H^{k+j}} \\ 
&\qquad \qquad + (1+t)^{-\frac{n}{2(1-\delta)}(\frac{1}{m}-\frac{1}{2})- \frac{k-2\delta}{2(1-\delta)}-j}\|w_1\|_{L^m \cap H^{[k+(j-1)]^+}}
\end{align*}
for all space dimensions $n>2m_0 \delta$.
\emd

\begin{nx}
\fontshape{n}
\selectfont
We recognize that the decay rates from Propositions \ref{md2.1} and \ref{md2.2} coincide with those in \cite{DabbiccoReissig}. Moreover, the optimality of those from Proposition \ref{md2.1} is also guaranteed by the study of asymptotic profile of solutions to (\ref{pt1.3}) in \cite{Dao2}. From this observation, these estimates play really a fundamental role in the proofs of global (in time) existence results for (\ref{pt1.1}) in Section \ref{Sec.Global existence}.
\end{nx}

\noindent Finally, plugging $\delta= \delta_\ell$ with $\ell=1,2$ into the statements from Propositions \ref{md2.1} and \ref{md2.2} we may arrive at the following result.
\bhq \label{hq2.1}
Let $\delta= \delta_\ell \in \big[0,\f{1}{2}\big]$ with $\ell=1,\,2$ in (\ref{pt1.3}). Let $k\ge 0$, $j=0,1$ and $m \in [1,2)$. The solutions to (\ref{pt1.3}) satisfy the following $(L^m \cap L^2)-L^2$ estimates:
\begin{align*}
\big\|\partial_t^j \nabla^k w(t,\cdot)\big\|_{L^2} &\lesssim (1+t)^{-\frac{n}{2(1-\delta_\ell)}(\frac{1}{m}-\frac{1}{2})- \frac{k}{2(1-\delta_\ell)}-j}\|w_0\|_{L^m \cap H^{k+j}} \\ 
&\qquad \qquad + (1+t)^{-\frac{n}{2(1-\delta_\ell)}(\frac{1}{m}-\frac{1}{2})- \frac{k-2\delta_\ell}{2(1-\delta_\ell)}-j}\|w_1\|_{L^m \cap H^{[k+(j-1)]^+}}
\end{align*}
for all space dimensions $n>2m_0 \max\{\delta_1,\, \delta_2\}$. Moreover, the following $L^2-L^2$ estimates hold:
$$ \big\|\partial_t^j \nabla^k w(t,\cdot)\big\|_{L^2} \lesssim (1+t)^{-\frac{k+2j\delta_\ell}{2(1-\delta_\ell)}}\|w_0\|_{H^{k+j}}+ (1+t)^{-\frac{k}{2(1-\delta_\ell)}-j+1}\|w_1\|_{H^{[k+(j-1)]^+}} $$
for all space dimensions $n \ge 1$.
\ehq

\subsection{A modified test function}
Main aim of this section is to provide some auxiliary properties of the modified test function $\phi= \phi(x):= \big< x\big>^{-r}$ for some $r>0$ from the recent paper of Dao-Reissig \cite{DaoReissig} which are key tools in the proof of our blow-up result in Section \ref{Sec.Blow-up}.
\bdn[\cite{Kwanicki,Silvestre}] \label{def1}
\fontshape{n}
\selectfont
Let $s \in (0,1)$. Let $X$ be a suitable set of functions defined on $\R^n$. Then, the fractional Laplacian $(-\Delta)^s$ in $\R^n$ is a non-local operator given by
$$ (-\Delta)^s: \,\,\phi \in X  \to (-\Delta)^s \phi(x):= C_{n,s}\,\, p.v.\int_{\R^n}\frac{\phi(x)- \phi(y)}{|x-y|^{n+2s}}dy $$
as long as the right-hand side exists, where $p.v.$ stands for Cauchy's principal value, $C_{n,s}:= \frac{4^s \Gamma(\frac{n}{2}+s)}{\pi^{\frac{n}{2}}\Gamma(-s)}$ is a normalization constant and $\Gamma$ denotes the Gamma function.
\edn

\bbd[Lemma 2.3 in \cite{DaoReissig} with $m=0$] \label{lemma2.1}
Let $s \in (0,1)$ and $r>0$. Then, the following estimates hold for all $x \in \R^n$:
$$ \big|(-\Delta)^s \big< x\big>^{-r}\big| \lesssim
\begin{cases}
\big< x\big>^{-r-2s} &\quad \text{ if }\quad 0< r< n, \\
\big< x\big>^{-n-2s}\log(e+|x|) &\quad \text{ if }\quad r= n, \\
\big< x\big>^{-n-2s} &\quad \text{ if }\quad r> n.
\end{cases} $$
\ebd

\bbd \label{lemma2.2}
Let $s\in (0,1)$. Let $\psi$ be a smooth function satisfying $\partial_x^2 \phi \in L^\ity$. For any $R>0$, let $\phi_R$ be a function defined by
$$ \phi_R(x):= \phi\big(R^{-\kappa} x\big)\quad \text{ for all }x \in \R^n, $$
where $\kappa>0$. Then, $(-\Delta)^s (\phi_R)$ satisfies the following scaling properties for all $x \in \R^n$:
$$(-\Delta)^s (\phi_R)(x)= R^{-2\kappa s}\big((-\Delta)^s \phi \big)\big(R^{-\kappa} x\big). $$
\ebd
\begin{proof}
We follow the proof of Lemma 2.4 in \cite{DaoReissig} with minor modifications to conclude the desired statement.
\end{proof}

\bbd[Lemma 2.7 in \cite{DaoReissig}] \label{lemma2.3}
Let $s \in \R$. Let $\phi_1= \phi_1(x) \in H^s$ and $\phi_2= \phi_2(x) \in H^{-s}$. Then, the following relation holds:
$$ \int_{\R^n}\phi_1(x)\,\phi_2(x)dx= \int_{\R^n}\widehat{\phi}_1(\xi)\,\widehat{\phi}_2(\xi)d\xi. $$
\ebd

\section{Global (in time) existence of small data solutions} \label{Sec.Global existence}

\subsection{Proof of Theorem \ref{dl1.1}}
At first, let recall the fundamental solutions
$$K_{0,\delta}(t,x)= \mathfrak{F}^{-1}_{\xi \to x}\big(\widehat{K}_{0,\delta}(t,\xi)\big) \text{ and } K_{1,\delta}(t,x)= \mathfrak{F}^{-1}_{\xi \to x}\big(\widehat{K}_{1,\delta}(t,\xi)\big) $$
defined in Section \ref{Sec.Pre} to represent the solutions of the corresponding linear Cauchy problems with vanishing right-hand sides to (\ref{pt1.1}) in the form
$$\begin{cases}
u^{ln}(t,x)=K_{0,\delta_1}(t,x) \ast_{x} u_0(x)+ K_{1,\delta_1}(t,x) \ast_{x} u_1(x), \\ 
v^{ln}(t,x)=K_{0,\delta_2}(t,x) \ast_{x} v_0(x)+ K_{1,\delta_2}(t,x) \ast_{x} v_1(x).
\end{cases}$$
By applying Duhamel's principle, the formal implicit representation of the solutions to (\ref{pt1.1}) can be read as follows:
$$\begin{cases}
u(t,x)= u^{ln}(t,x) + \displaystyle\int_0^t K_{1,\delta_1}(t-\tau,x) \ast_x |v(\tau,x)|^p d\tau=: u^{ln}(t,x)+ u^{nl}(t,x), \\ 
v(t,x)= v^{ln}(t,x) + \displaystyle\int_0^t K_{1,\delta_2}(t-\tau,x) \ast_x |u(\tau,x)|^q d\tau=: v^{ln}(t,x)+ v^{nl}(t,x).
\end{cases}$$
Let us now choose the data spaces $(u_0,u_1) \in \mathcal{A}$ and $(v_0,v_1) \in \mathcal{A}$. We introduce the family $\{X(t)\}_{t>0}$ of the solution spaces
$$ X(t):= \Big(C\big([0,t],H^1\big)\cap C^1\big([0,t],L^2\big)\Big)^2, $$
equipped with the following norm:
\begin{align*}
\|(u,v)\|_{X(t)}:= \sup_{0\le \tau \le t} \Big( &f_1(\tau)^{-1}\|u(\tau,\cdot)\|_{L^2} + f_2(\tau)^{-1}\big\|\nabla u(\tau,\cdot)\big\|_{L^2}+ f_3(\tau)^{-1}\|u_t(\tau,\cdot)\|_{L^2} \\
&+g_1(\tau)^{-1}\|v(\tau,\cdot)\|_{L^2} + g_2(\tau)^{-1}\big\|\nabla v(\tau,\cdot)\big\|_{L^2}+ g_3(\tau)^{-1}\|v_t(\tau,\cdot)\|_{L^2} \Big),
\end{align*}
where
\begin{align}
&f_1(\tau)= (1+\tau)^{-\frac{n}{2(1-\delta_1)}(\frac{1}{m}-\frac{1}{2})+ \frac{\delta_1}{1-\delta_1}+ \e(p,\delta_2)},\quad f_2(\tau)= (1+\tau)^{-\frac{n}{2(1-\delta_1)}(\frac{1}{m}-\frac{1}{2})- \frac{1-2\delta_1}{2(1-\delta_1)}+ \e(p,\delta_2)}, \label{pt4.11} \\
&f_3(\tau)=(1+\tau)^{-\frac{n}{2(1-\delta_1)}(\frac{1}{m}-\frac{1}{2})- \frac{1-2\delta_1}{1-\delta_1}+ \e(p,\delta_2)}, \label{pt4.12} \\
&g_1(\tau)= (1+\tau)^{-\frac{n}{2(1-\delta_2)}(\frac{1}{m}-\frac{1}{2})+ \frac{\delta_2}{1-\delta_2}},\quad g_2(\tau)= (1+\tau)^{-\frac{n}{2(1-\delta_2)}(\frac{1}{m}-\frac{1}{2})- \frac{1-2\delta_2}{2(1-\delta_2)}}, \label{pt4.21} \\
&g_3(\tau)=(1+\tau)^{-\frac{n}{2(1-\delta_2)}(\frac{1}{m}-\frac{1}{2})- \frac{1-2\delta_2}{1-\delta_2}}. \label{pt4.22}
\end{align}
For all $t>0$, we define the following operator:
\begin{align*}
N &: \quad X(t) \longrightarrow X(t) \\
N(u,v)(t,x) &= \big(u^{ln}(t,x)+ u^{nl}(t,x), v^{ln}(t,x)+ v^{nl}(t,x)\big).
\end{align*}
Our main goal is to prove the operator $N$ satisfying the following two inequalities:
\begin{align}
\|N(u,v)\|_{X(t)} &\lesssim \|(u_0,u_1)\|_{\mathcal{A}}+ \|(v_0,v_1)\|_{\mathcal{A}}+ \|(u,v)\|^{p}_{X(t)}+ \|(u,v)\|^{q}_{X(t)}, \label{pt4.3} \\
\|N(u,v)-N(\bar{u},\bar{v})\|_{X(t)} &\lesssim \|(u,v)-(\bar{u},\bar{v})\|_{X(t)} \Big(\|(u,v)\|^{p-1}_{X(t)}+ \|(\bar{u},\bar{v})\|^{p-1}_{X(t)} \nonumber \\
&\hspace{6cm}+ \|(u,v)\|^{q-1}_{X(t)}+ \|(\bar{u},\bar{v})\|^{q-1}_{X(t)}\Big), \label{pt4.4}
\end{align}
Then, employing Banach's fixed point theorem we may conclude global (in time) existence results of small data solutions. For this purpose, we replace $j,k=0,1$ with $(j,k)\ne (1,1)$ into Corollary \ref{hq2.1} to arrive at the estimate
\begin{equation*}
\big\|(u^{ln}, v^{ln})\big\|_{X(t)} \lesssim \|(u_0,u_1)\|_{\mathcal{A}}+ \|(v_0,v_1)\|_{\mathcal{A}}
\end{equation*}
by the definition of the norm in $X(t)$. Hence, it is reasonable to prove the following inequality instead of (\ref{pt4.3}):
\begin{equation} \label{pt4.5}
\big\|(u^{nl}, v^{nl})\big\|_{X(t)} \lesssim \|(u,v)\|^{p}_{X(t)}+ \|(u,v)\|^{q}_{X(t)}.
\end{equation}
\par First, let us show the inequality (\ref{pt4.5}). To deal with $u^{nl}$, we use the $(L^m \cap L^2)- L^2$ estimates from Corollary \ref{hq2.1} to get the following estimate:
$$ \big\|u^{nl}(t,\cdot)\big\|_{L^2}\lesssim \int_0^t (1+t-\tau)^{-\frac{n}{2(1-\delta_1)}(\frac{1}{m}-\frac{1}{2})+ \frac{\delta_1}{1-\delta_1}}\big\||v(\tau,\cdot)|^p\big\|_{L^m \cap L^2}d\tau. $$
Thus, we need to estimate $|v(\tau,x)|^p$ in $L^m$ and $L^2$. We have
$$ \big\||v(\tau,\cdot)|^p\big\|_{L^m}= \|v(\tau,\cdot)\|^p_{L^{mp}} \quad \text{ and }\quad \big\||v(\tau,\cdot)|^p\big\|_{L^2}= \|v(\tau,\cdot)\|^p_{L^{2p}}. $$
After using the fractional Gagliardo-Nirenberg inequality from Proposition \ref{fractionalGagliardoNirenberg}, we deduce that
\begin{align}
\big\||v(\tau,\cdot)|^p\big\|_{L^m} &\lesssim (1+\tau)^{-\frac{n}{2m(1-\delta_2)}(p-1)+ \frac{p\delta_2}{1-\delta_2}}\|(u,v)\|^p_{X(\tau)}, \label{t11.1} \\
\big\||v(\tau,\cdot)|^p\big\|_{L^2} &\lesssim (1+\tau)^{-\frac{n}{2m(1-\delta_2)}(p-\frac{1}{2})+ \frac{p\delta_2}{1-\delta_2}}\|(u,v)\|^p_{X(\tau)}, \label{t11.2}
\end{align}
where the conditions (\ref{GN11A1}) and (\ref{GN11A2}) are fulfilled for $p$. As a consequence, we can proceed as follows:
\begin{align*}
\big\|u^{nl}(t,\cdot)\big\|_{L^2} &\lesssim (1+t)^{-\frac{n}{2(1-\delta_1)}(\frac{1}{m}-\frac{1}{2})+ \frac{\delta_1}{1-\delta_1}}\|(u,v)\|^p_{X(t)} \int_0^{t/2}(1+\tau)^{-\frac{n}{2m(1-\delta_2)}(p-1)+ \frac{p\delta_2}{1-\delta_2}} d\tau \\
&\qquad + (1+t)^{-\frac{n}{2m(1-\delta_2)}(p-1)+ \frac{p\delta_2}{1-\delta_2}}\|(u,v)\|^p_{X(t)} \int_{t/2}^t (1+t-\tau)^{-\frac{n}{2(1-\delta_1)}(\frac{1}{m}-\frac{1}{2})+ \frac{\delta_1}{1-\delta_1}}d\tau,
\end{align*}
where we notice that we used the relation
\begin{equation}
(1+t-\tau) \approx (1+t) \text{ if }\tau \in [0,t/2], \text{ and } (1+\tau) \approx (1+t) \text{ if }\tau \in [t/2,t]. \label{t11.3}
\end{equation}
Due to the condition $p \le 1+ \frac{2m}{n- 2m\delta_2}$ in (\ref{exponent11A2}), it implies immediately that the term $(1+\tau)^{-\frac{n}{2m(1-\delta_2)}(p-1)+ \frac{p\delta_2}{1-\delta_2}}$ is not integrable. For this reason, may estimate
\begin{align*}
&(1+t)^{-\frac{n}{2(1-\delta_1)}(\frac{1}{m}-\frac{1}{2})+ \frac{\delta_1}{1-\delta_1}} \int_0^{t/2}(1+\tau)^{-\frac{n}{2m(1-\delta_2)}(p-1)+ \frac{p\delta_2}{1-\delta_2}} d\tau \\
&\qquad \lesssim
\begin{cases}
(1+t)^{-\frac{n}{2(1-\delta_1)}(\frac{1}{m}-\frac{1}{2})+ \frac{\delta_1}{1-\delta_1}+1-\frac{n}{2m(1-\delta_2)}(p-1)+ \frac{p\delta_2}{1-\delta_2}} &\text{ if }p< 1+ \frac{2m}{n- 2m\delta_2} \\
(1+t)^{-\frac{n}{2(1-\delta_1)}(\frac{1}{m}-\frac{1}{2})+ \frac{\delta_1}{1-\delta_1}+\e} &\text{ if }p= 1+ \frac{2m}{n- 2m\delta_2}
\end{cases} \\ 
&\qquad \lesssim (1+t)^{-\frac{n}{2(1-\delta_1)}(\frac{1}{m}-\frac{1}{2})+ \frac{\delta_1}{1-\delta_1}+ \e(p,\delta_2)},
\end{align*}
where $\e$ is a sufficiently small positive number. Thanks to the condition $n\le \frac{4m}{2-m}$ in (\ref{GN11A2}), we may verify that $-\frac{n}{2(1-\delta_1)}(\frac{1}{m}-\frac{1}{2})+ \frac{\delta_1}{1-\delta_1}\ge -1$. Hence, we derive
\begin{align*}
&(1+t)^{-\frac{n}{2m(1-\delta_2)}(p-1)+ \frac{p\delta_2}{1-\delta_2}} \int_{t/2}^t (1+t-\tau)^{-\frac{n}{2(1-\delta_1)}(\frac{1}{m}-\frac{1}{2})+ \frac{\delta_1}{1-\delta_1}}d\tau \\
&\qquad \lesssim
\begin{cases}
(1+t)^{-\frac{n}{2m(1-\delta_2)}(p-1)+ \frac{p\delta_2}{1-\delta_2}+1-\frac{n}{2(1-\delta_1)}(\frac{1}{m}-\frac{1}{2})+ \frac{\delta_1}{1-\delta_1}} &\text{ if } n< \frac{4m}{2-m} \\
(1+t)^{-\frac{n}{2m(1-\delta_2)}(p-1)+ \frac{p\delta_2}{1-\delta_2}+\e} &\text{ if } n= \frac{4m}{2-m}
\end{cases} \\
&\qquad \lesssim (1+t)^{-\frac{n}{2(1-\delta_1)}(\frac{1}{m}-\frac{1}{2})+ \frac{\delta_1}{1-\delta_1}+ \e(p,\delta_2)},
\end{align*}
where  $\e$ is a sufficiently small positive number. Therefore, combining the above estimates we may conclude the following estimate:
$$\big\|u^{nl}(t,\cdot)\big\|_{L^2} \lesssim (1+t)^{-\frac{n}{2(1-\delta_1)}(\frac{1}{m}-\frac{1}{2})+ \frac{\delta_1}{1-\delta_1}+ \e(p,\delta_2)} \|(u,v)\|^p_{X(t)}. $$
In order to control $\nabla u^{nl}$, we use the $(L^m \cap L^2)- L^2$ estimates if $\tau \in [0,t/2]$ and the $L^2-L^2$ estimates if $\tau \in [t/2,t]$ from Corollary \ref{hq2.1} to arrive at
\begin{align*}
\big\|\nabla u^{nl}(t,\cdot)\big\|_{L^2}&\lesssim \int_0^{t/2} (1+t-\tau)^{-\frac{n}{2(1-\delta_1)}(\frac{1}{m}-\frac{1}{2})- \frac{1-2\delta_1}{2(1-\delta_1)}}\big\||v(\tau,\cdot)|^p\big\|_{L^m \cap L^2}d\tau \\
&\qquad +\int_{t/2}^t (1+t-\tau)^{-\frac{1-2\delta_1}{2(1-\delta_1)}}\big\||v(\tau,\cdot)|^p\big\|_{L^2}d\tau \\
&\lesssim (1+t)^{-\frac{n}{2(1-\delta_1)}(\frac{1}{m}-\frac{1}{2})- \frac{1-2\delta_1}{2(1-\delta_1)}}\|(u,v)\|^p_{X(t)} \int_0^{t/2}(1+\tau)^{-\frac{n}{2m(1-\delta_2)}(p-1)+ \frac{p\delta_2}{1-\delta_2}}d\tau \\
&\qquad +(1+t)^{-\frac{n}{2m(1-\delta_2)}(p-\frac{1}{2})+ \frac{p\delta_2}{1-\delta_2}}\|(u,v)\|^p_{X(t)} \int_{t/2}^t (1+t-\tau)^{-\frac{1-2\delta_1}{2(1-\delta_1)}}d\tau,
\end{align*}
where we used again the estimates (\ref{t11.1}) and (\ref{t11.2}) linked to the relation (\ref{t11.3}). In the same treatment of $u^{nl}$, we obtain the following estimate for first integral:
\begin{align*}
&(1+t)^{-\frac{n}{2(1-\delta_1)}(\frac{1}{m}-\frac{1}{2})- \frac{1-2\delta_1}{2(1-\delta_1)}} \int_0^{t/2}(1+\tau)^{-\frac{n}{2m(1-\delta_2)}(p-1)+ \frac{p\delta_2}{1-\delta_2}}d\tau \\ 
&\qquad \lesssim (1+t)^{-\frac{n}{2(1-\delta_1)}(\frac{1}{m}-\frac{1}{2})- \frac{1-2\delta_1}{2(1-\delta_1)}+ \e(p,\delta_2)}.
\end{align*}
Moreover, the remaining integral can be dealt with the following way:
\begin{align*}
&(1+t)^{-\frac{n}{2m(1-\delta_2)}(p-\frac{1}{2})+ \frac{p\delta_2}{1-\delta_2}} \int_{t/2}^t (1+t-\tau)^{-\frac{1-2\delta_1}{2(1-\delta_1)}}d\tau \\ 
&\qquad \lesssim (1+t)^{-\frac{n}{2m(1-\delta_2)}(p-\frac{1}{2})+ \frac{p\delta_2}{1-\delta_2}+1 -\frac{1-2\delta_1}{2(1-\delta_1)}}\lesssim (1+t)^{-\frac{n}{2(1-\delta_1)}(\frac{1}{m}-\frac{1}{2})- \frac{1-2\delta_1}{2(1-\delta_1)}+ \e(p,\delta_2)}
\end{align*}
due to $\delta_1 \ge \delta_2$. Consequently, we have shown that
$$\big\|\nabla u^{nl}(t,\cdot)\big\|_{L^2} \lesssim (1+t)^{-\frac{n}{2(1-\delta_1)}(\frac{1}{m}-\frac{1}{2})- \frac{1-2\delta_1}{2(1-\delta_1)}+ \e(p,\delta_2)} \|(u,v)\|^p_{X(t)}. $$
By analogous arguments as we estimated $\nabla u^{nl}$ we also derive
$$\big\|u_t^{nl}(t,\cdot)\big\|_{L^2} \lesssim (1+t)^{-\frac{n}{2(1-\delta_1)}(\frac{1}{m}-\frac{1}{2})- \frac{1-2\delta_1}{1-\delta_1}+ \e(p,\delta_2)}\|(u,v)\|^p_{X(t)}. $$
Similarly, we may conclude the following estimates for $j,k=0,1$ with $(j,k)\ne (1,1)$:
$$\big\|\partial^j_t \nabla^k v^{nl}(t,\cdot)\big\|_{L^2} \lesssim (1+t)^{-\frac{n}{2(1-\delta_2)}(\frac{1}{m}-\frac{1}{2})- \frac{k-2\delta_2}{2(1-\delta_2)}-j}\|(u,v)\|^q_{X(t)}, $$
provided that the conditions from (\ref{GN11A1}) to (\ref{exponent11A2}) are satisfied for $q$. Therefore, from the definition of the norm in $X(t)$ we have proved that the inequality (\ref{pt4.5}) holds.
\par Let us now indicate the inequality (\ref{pt4.4}). For two elements $(u,v)$ and $(\bar{u},\bar{v})$ from $X(t)$, we get
$$N(u,v)(t,x)- N(\bar{u},\bar{v})(t,x)= \big(u^{nl}(t,x)- \bar{u}^{nl}(t,x), v^{nl}(t,x)- \bar{v}^{nl}(t,x)\big). $$
The proof of (\ref{pt4.4}) can be proceeded in the same ways as that of (\ref{pt4.5}). For this reason, let us sketch our proof. On the one hand, we use the $(L^m \cap L^2)- L^2$ estimates from Corollary \ref{hq2.1} for $u^{nl}- \bar{u}^{nl}$ and $v^{nl}- \bar{v}^{nl}$. Meanwhile, for $\partial^j_t \nabla^k \big(u^{nl}- \bar{u}^{nl}\big)$ and $\partial^j_t \nabla^k \big(v^{nl}- \bar{v}^{nl}\big)$, with $(j,k)=(0,1)$ or $(1,0)$, we apply $(L^m \cap L^2)- L^2$ estimates if $\tau \in [0,t/2]$ and the $L^2-L^2$ estimates if $\tau \in [t/2,t]$ from Corollary \ref{hq2.1}. Therefore, we may arrive at the following estimates for $(j,k)=(0,1)$ or $(1,0)$:
\begin{align*}
\big\|\big(u^{nl}- \bar{u}^{nl}\big)(t,\cdot)\big\|_{L^2} &\lesssim \int_0^t (1+t-\tau)^{-\frac{n}{2(1-\delta_1)}(\frac{1}{m}-\frac{1}{2})+ \frac{\delta_1}{1-\delta_1}}\big\||v(\tau,\cdot)|^p- \bar{v}(\tau,\cdot)|^p\big\|_{L^m \cap L^2}d\tau, \\
\big\|\partial^j_t \nabla^k \big(u^{nl}- \bar{u}^{nl}\big)(t,\cdot)\big\|_{L^2} &\lesssim \int_0^{t/2}(1+t-\tau)^{-\frac{n}{2(1-\delta_1)}(\frac{1}{m}-\frac{1}{2})- \frac{k-2\delta_1}{2(1-\delta_1)}-j}\big\||v(\tau,\cdot)|^p- \bar{v}(\tau,\cdot)|^p\big\|_{L^m \cap L^2}d\tau \\
&\qquad + \int_{t/2}^t (1+t-\tau)^{-\frac{k-2\delta_1}{2(1-\delta_1)}-j}\big\||v(\tau,\cdot)|^p- \bar{v}(\tau,\cdot)|^p\big\|_{L^2}d\tau,
\end{align*}
and
\begin{align*}
\big\|\big(v^{nl}- \bar{v}^{nl}\big)(t,\cdot)\big\|_{L^2} &\lesssim \int_0^t (1+t-\tau)^{-\frac{n}{2(1-\delta_2)}(\frac{1}{m}-\frac{1}{2})+ \frac{\delta_2}{1-\delta_2}}\big\||u(\tau,\cdot)|^q- \bar{u}(\tau,\cdot)|^q\big\|_{L^m \cap L^2}d\tau, \\
\big\|\partial^j_t \nabla^k \big(v^{nl}- \bar{v}^{nl}\big)(t,\cdot)\big\|_{L^2} &\lesssim \int_0^{t/2}(1+t-\tau)^{-\frac{n}{2(1-\delta_2)}(\frac{1}{m}-\frac{1}{2})- \frac{k-2\delta_2}{2(1-\delta_2)}-j}\big\||u(\tau,\cdot)|^q- \bar{u}(\tau,\cdot)|^q\big\|_{L^m \cap L^2}d\tau, \\
&\qquad + \int_{t/2}^t (1+t-\tau)^{-\frac{k-2\delta_2}{2(1-\delta_2)}-j}\big\||u(\tau,\cdot)|^q- \bar{u}(\tau,\cdot)|^q\big\|_{L^2}d\tau.
\end{align*}
By applying H\"{o}lder's inequality we have
\begin{align*}
\big\||v(\tau,\cdot)|^p- |\bar{v}(\tau,\cdot)|^p\big\|_{L^m} &\lesssim \|v(\tau,\cdot)- \bar{v}(\tau,\cdot)\|_{L^{mp}} \big(\|v(\tau,\cdot)\|^{p-1}_{L^{mp}}+ \|\bar{v}(\tau,\cdot)\|^{p-1}_{L^{mp}}\big), \\
\big\||v(\tau,\cdot)|^p- |\bar{v}(\tau,\cdot)|^p\big\|_{L^2} &\lesssim \|v(\tau,\cdot)- \bar{v}(\tau,\cdot)\|_{L^{2p}} \big(\|v(\tau,\cdot)\|^{p-1}_{L^{2p}}+ \|\bar{v}(\tau,\cdot)\|^{p-1}_{L^{2p}}\big), \\
\big\||u(\tau,\cdot)|^q- |\bar{u}(\tau,\cdot)|^q\big\|_{L^m} &\lesssim \|u(\tau,\cdot)- \bar{u}(\tau,\cdot)\|_{L^{mq}} \big(\|u(\tau,\cdot)\|^{q-1}_{L^{mq}}+ \|\bar{u}(\tau,\cdot)\|^{q-1}_{L^{mq}}\big), \\
\big\||u(\tau,\cdot)|^q- |\bar{u}(\tau,\cdot)|^q\big\|_{L^2} &\lesssim \|u(\tau,\cdot)- \bar{u}(\tau,\cdot)\|_{L^{2q}} \big(\|u(\tau,\cdot)\|^{q-1}_{L^{2q}}+ \|\bar{u}(\tau,\cdot)\|^{q-1}_{L^{2q}}\big).
\end{align*}
Analogously to the proof of (\ref{pt4.5}), we employ the fractional Gagliardo-Nirenberg inequality from Proposition \ref{fractionalGagliardoNirenberg} to the terms
\begin{align*}
&\|v(\tau,\cdot)-\bar{v}(\tau,\cdot)\|_{L^{\eta_1}},\quad \|u(\tau,\cdot)-\bar{u}(\tau,\cdot)\|_{L^{\eta_2}}, \\
&\|v(\tau,\cdot)\|_{L^{\eta_1}},\quad \|\bar{v}(\tau,\cdot)\|_{L^{\eta_1}},\quad \|u(\tau,\cdot)\|_{L^{\eta_2}},\quad \|\bar{u}(\tau,\cdot)\|_{L^{\eta_2}},
\end{align*}
with $\eta_1= mp$ or $\eta_1= 2p$, and $\eta_2= mq$ or $\eta_2=2q$ to complete the proof of the inequality (\ref{pt4.4}). Summarizing, Theorem \ref{dl1.1} is proved completedly.

\subsection{Proof of Theorem \ref{dl1.2}}
We follow the proof of Theorem \ref{dl1.1} with minor modifications in the steps of our proof. We also introduce both spaces for the data and the solutions as in Theorem \ref{dl1.1}, where the weights (\ref{pt4.11}) to (\ref{pt4.22}) are modified in the following way:
\begin{align*}
&f_1(\tau)= (1+\tau)^{-\frac{n}{2(1-\delta_1)}(\frac{1}{m}-\frac{1}{2})+ \frac{\delta_1}{1-\delta_1}},\quad f_2(\tau)= (1+\tau)^{-\frac{n}{2(1-\delta_1)}(\frac{1}{m}-\frac{1}{2})- \frac{1-2\delta_1}{2(1-\delta_1)}}, \\
&f_3(\tau)=(1+\tau)^{-\frac{n}{2(1-\delta_1)}(\frac{1}{m}-\frac{1}{2})- \frac{1-2\delta_1}{1-\delta_1}}, \\
&g_1(\tau)= (1+\tau)^{-\frac{n}{2(1-\delta_2)}(\frac{1}{m}-\frac{1}{2})+ \frac{\delta_2}{1-\delta_2}+ \e(q,\delta_1)},\quad g_2(\tau)= (1+\tau)^{-\frac{n}{2(1-\delta_2)}(\frac{1}{m}-\frac{1}{2})- \frac{1-2\delta_2}{2(1-\delta_2)}+ \e(q,\delta_1)}, \\
&g_3(\tau)=(1+\tau)^{-\frac{n}{2(1-\delta_2)}(\frac{1}{m}-\frac{1}{2})- \frac{1-2\delta_2}{1-\delta_2}+ \e(q,\delta_1)}.
\end{align*}
Then, repeating some steps of the proofs we did in Theorem \ref{dl1.1} we may conclude the proof of Theorem \ref{dl1.2}.

\section{Nonexistence result via modified test function method} \label{Sec.Blow-up}
In order to prove the blow-up results, we shall apply a modified test function method from Section \ref{Sec.Pre} which plays a significant role in the following proofs.

\subsection{Proof of Theorem \ref{dl1.3}}
First, we introduce the function $\varphi= \varphi(t)$ having the following properties:
\begin{align}
&1.\quad \varphi \in \mathcal{C}_0^\ity([0,\ity)) \text{ and }
\varphi(t)=\begin{cases}
1 &\quad \text{ for }0 \le t \le \frac{1}{2}, \\
\text{decreasing } &\quad \text{ for }\frac{1}{2} \le t \le 1, \\
0 &\quad \text{ for }t \ge 1,
\end{cases} \nonumber \\
&2.\quad \varphi^{-\frac{\kappa'}{\kappa}}(t)\big(|\varphi'(t)|^{\kappa'}+|\varphi''(t)|^{\kappa'}\big) \le C \quad \text{ for any } t \in \Big[\frac{1}{2},1\Big], \label{t13.1}
\end{align}
with $\kappa= p$ or $\kappa= q$, where $\kappa'$ is the conjugate of $\kappa>1$ and $C$ is a suitable positive constant. Now we denote $\delta_0:= \min\{\delta_1,\,\delta_2\}$. Due to the assumption of both $\delta_1$ and $\delta_2 \in (0,1)$, it is clear that $\delta_0 \in (0,1)$, too. Then, we introduce the function $\psi= \psi(|x|):= \big< x\big>^{-n-2\delta_0}$.
\par Let $R$ be a large parameter in $[0,\ity)$. We define the following test function:
$$ \eta_R(t,x):= \varphi_R(t) \psi_R(x), $$
where $\varphi_R(t):= \varphi(R^{-\alpha}t)$ and $\psi_R(x):= \psi(R^{-\beta}x)$ for some $\alpha,\,\beta$ which we will fix later. We define the functionals
\begin{align*}
&I_R:= \int_0^{\ity}\int_{\R^n}|v(t,x)|^p \eta_R(t,x)\,dxdt= \int_0^{R^{\alpha}}\int_{\R^n}|v(t,x)|^p \eta_R(t,x)\,dxdt, \\
&J_R:=  \int_0^{\ity}\int_{\R^n}|u(t,x)|^q \eta_R(t,x)\,dxdt= \int_0^{R^{\alpha}}\int_{\R^n}|u(t,x)|^q \eta_R(t,x)\,dxdt,
\end{align*}
and
$$ I_{R,t}:= \int_{\frac{R^\alpha}{2}}^{R^{\alpha}}\int_{\R^n}|v(t,x)|^p \eta_R(t,x)\,dxdt, \quad J_{R,t}:= \int_{\frac{R^\alpha}{2}}^{R^{\alpha}}\int_{\R^n}|u(t,x)|^q \eta_R(t,x)\,dxdt. $$
Let us assume that $(u,v)= \big(u(t,x),v(t,x)\big)$ is a global (in time) Sobolev solution from $C\big([0,\infty),L^2\big) \times C\big([0,\infty),L^2\big)$ to (\ref{pt1.1}). We multiply the first equation to (\ref{pt1.1}) by $\eta_R=\eta_R(t,x)$ and carry out partial integration to get
\begin{align}
0\le I_R &= -\int_{\R^n} u_1(x)\psi_R(x)\,dx + \int_{\frac{R^\alpha}{2}}^{R^{\alpha}}\int_{\R^n}u(t,x) \varphi''_R(t) \psi_R(x)\,dxdt \nonumber \\
&\quad - \int_0^{\ity}\int_{\R^n} \varphi_R(t) \psi_R(x)\, \Delta u(t,x)\,dxdt- \int_{\frac{R^\alpha}{2}}^{R^\alpha}\int_{\R^n} \varphi'_R(t) \psi_R(x)\,(-\Delta)^{\delta_1} u(t,x)\,dxdt \nonumber \\
&=: -\int_{\R^n} u_1(x)\psi_R(x)\,dx+ I_{1R}- I_{2R}- I_{3R}. \label{t13.2}
\end{align}
Employing H\"{o}lder's inequality with $\frac{1}{q}+\frac{1}{q'}=1$ we can proceed as follows:
\begin{align*}
|I_{1R}| &\le \int_{\frac{R^\alpha}{2}}^{R^{\alpha}}\int_{\R^n} |u(t,x)|\, \big|\varphi''_R(t)\big| \psi_R(x) \, dxdt \\
&\lesssim \Big(\int_{\frac{R^\alpha}{2}}^{R^{\alpha}}\int_{\R^n} \Big|u(t,x)\eta^{\frac{1}{q}}_R(t,x)\Big|^p \,dxdt\Big)^{\frac{1}{q}} \Big(\int_{\frac{R^\alpha}{2}}^{R^{\alpha}}\int_{\R^n} \Big|\eta^{-\frac{1}{q}}_R(t,x) \varphi''_R(t) \psi_R(x)\Big|^{p'}\, dxdt\Big)^{\frac{1}{q'}} \\
&\lesssim J_{R,t}^{\frac{1}{q}}\, \Big(\int_{\frac{R^\alpha}{2}}^{R^{\alpha}}\int_{\R^n} \varphi_R^{-\frac{q'}{q}}(t) \big|\varphi''_R(t)\big|^{q'} \psi_R(x)\, dxdt\Big)^{\frac{1}{q'}}.
\end{align*}
After performing the change of variables $\tilde{t}:= R^{-\alpha}t$ and $\tilde{x}:= R^{-\beta}x$, we calculate straightforwardly to obtain
\begin{equation}
|I_{1R}| \lesssim J_{R,t}^{\frac{1}{q}}\, R^{-2\alpha+ \frac{\alpha+n\beta}{q'}}\Big(\int_{\R^n} \big< \tilde{x}\big>^{-n-2\delta_0}\, d\tilde{x}\Big)^{\frac{1}{q'}}, \label{t13.3}
\end{equation}
where we used $\varphi''_R(t)= R^{-2\alpha}\varphi''(\tilde{t})$ and the assumption (\ref{t13.1}). Now let us focus our considerations to deal with $I_{2R}$ and $I_{3R}$. First, since $\psi_R \in H^2$ and $u \in C\big([0,\infty),L^2\big)$, we apply Lemma \ref{lemma2.3} to arrive at the following relations:
\begin{align*}
\int_{\R^n} \psi_R(x)\, (-\Delta) u(t,x)\,dx&= \int_{\R^n}|\xi|^2 \widehat{\psi}_R(\xi)\,\widehat{u}(t,\xi)\,d\xi= \int_{\R^n} u(t,x)\, (-\Delta) \psi_R(x)\,dx, \\
\int_{\R^n} \psi_R(x)\,(-\Delta)^{\delta_1} u(t,x)\,dx&= \int_{\R^n}|\xi|^{2\delta_1}\widehat{\psi}_R(\xi)\,\widehat{u}(t,\xi)\,d\xi= \int_{\R^n} u(t,x)\,(-\Delta)^{\delta_1}\psi_R(x)\,dx.
\end{align*}
As a consequence, it implies immediately that
$$I_{2R}= \int_0^{\ity}\int_{\R^n} \varphi_R(t) \psi_R(x)\, \Delta u(t,x)\,dxdt= \int_0^{\ity}\int_{\R^n} \varphi_R(t) u(t,x)\, \Delta \psi_R(x) \,dxdt, $$
and
$$I_{3R}= \int_{\frac{R^\alpha}{2}}^{R^\alpha}\int_{\R^n} \varphi'_R(t) \psi_R(x)\,(-\Delta)^{\delta_1} u(t,x)\,dxdt= \int_{\frac{R^\alpha}{2}}^{R^\alpha}\int_{\R^n} \varphi'_R(t) u(t,x)\,(-\Delta)^{\delta_1}\psi_R(x)\, dxdt. $$
Applying H\"{o}lder's inequality again as we estimated $J_1$ leads to
$$|I_{2R}|\le I_R^{\frac{1}{q}}\, \Big(\int_0^{R^{\alpha}}\int_{\R^n} \varphi_R(t) \psi^{-\frac{q'}{q}}_R(x)\, \big|\Delta \psi_R(x)\big|^{q'} \, dxdt\Big)^{\frac{1}{q'}}, $$
and
$$|I_{3R}|\le I_{R,t}^{\frac{1}{q}}\, \Big(\int_{\frac{R^\alpha}{2}}^{R^{\alpha}}\int_{\R^n} \varphi^{-\frac{q'}{q}}_R(t) \big|\varphi'_R(t)\big|^{q'} \psi^{\frac{-q'}{p}}_R(x)\, \big|(-\Delta)^{\delta_1}\psi_R(x)\big|^{q'} \, dxdt\Big)^{\frac{1}{q'}}. $$
To estimate the above two integrals, the key tools rely on results from Lemmas \ref{lemma2.1} and \ref{lemma2.2}. More in detail, in the first step we use the change of variables $\tilde{t}:= R^{-\alpha}t$ and $\tilde{x}:= R^{-\beta}x$ to derive
\begin{align*}
|I_{2R}| &\lesssim J_R^{\frac{1}{q}}\, R^{-2\beta+ \frac{\alpha+n\beta}{q'}}\Big(\int_0^{1}\int_{\R^n} \varphi(\tilde{t}) \psi^{-\frac{q'}{q}}(\tilde{x})\, \big|\Delta (\psi)(\tilde{x})\big|^{q'}\, d\tilde{x}d\tilde{t}\Big)^{\frac{1}{q'}} \\
&\lesssim J_R^{\frac{1}{q}}\, R^{-2\beta+ \frac{\alpha+n\beta}{q'}}\Big(\int_{\R^n} \psi^{-\frac{q'}{q}}(\tilde{x})\, \big|\Delta (\psi)(\tilde{x})\big|^{q'}\, d\tilde{x}\Big)^{\frac{1}{q'}},
\end{align*}
where we notice that $\Delta \psi_R(x)= R^{-2\beta} \Delta \varphi(\tilde{x})$. Hence, we deduce the following estimate:
\begin{equation}
|I_{2R}|\lesssim I_R^{\frac{1}{q}}\, R^{-2\beta+ \frac{\alpha+n\beta}{q'}}\Big(\int_{\R^n} \big< \tilde{x}\big>^{-n-2\delta_0-2q'}\, d\tilde{x}\Big)^{\frac{1}{q'}}. \label{t13.4}
\end{equation}
Now let us come back to estimate $I_{3R}$ in the second step. After carrying out again the change of variables $\tilde{t}:= R^{-\alpha}t$ and $\tilde{x}:= R^{-\beta}x$ and applying Lemma \ref{lemma2.2}, we may estimate $I_{3R}$ by
\begin{align*}
|I_{3R}| &\lesssim J_{R,t}^{\frac{1}{q}}\, R^{-\alpha-2\delta_1\beta+ \frac{\alpha+n\beta}{q'}}\Big(\int_{\frac{1}{2}}^{1}\int_{\R^n} \varphi^{-\frac{q'}{q}}(\tilde{t}) \big|\varphi'(\tilde{t})\big|^{q'} \psi^{-\frac{q'}{q}}(\tilde{x})\, \big|(-\Delta)^{\delta_1}(\psi)(\tilde{x})\big|^{q'}\, d\tilde{x}d\tilde{t}\Big)^{\frac{1}{q'}} \\
&\lesssim J_{R,t}^{\frac{1}{p}}\, R^{-\alpha-2\delta_1\beta+ \frac{\alpha+n\beta}{q'}}\Big(\int_{\R^n} \psi^{-\frac{q'}{q}}(\tilde{x})\, \big|(-\Delta)^{\delta_1}(\psi)(\tilde{x})\big|^{q'}\, d\tilde{x}\Big)^{\frac{1}{q'}},
\end{align*}
where we used $\varphi'_R(t)= R^{-\alpha}\varphi'(\tilde{t})$ and the assumption (\ref{t13.1}). In order to control the last integral, we employ Lemma \ref{lemma2.1} with $q=n+2\delta_0$ and $s=\delta_1$ to have
\begin{align}
|I_{3R}| &\lesssim J_{R,t}^{\frac{1}{q}} R^{-\alpha-2\delta_1\beta+ \frac{\alpha+n\beta}{q'}}\Big(\int_{\R^n} \big< \tilde{x}\big>^{-\frac{q'}{q}(-n-2\delta_0)}\, \big< \tilde{x}\big>^{q'(-n-2\delta_1)}\,d\tilde{x}\Big)^{\frac{1}{q'}} \nonumber \\ 
&\lesssim J_{R,t}^{\frac{1}{q}} R^{-\alpha-2\delta_1\beta+ \frac{\alpha+n\beta}{q'}}\Big(\int_{\R^n} \big< \tilde{x}\big>^{-n-2\delta_0}\,d\tilde{x}\Big)^{\frac{1}{q'}}. \label{t13.5}
\end{align}
Thanks to the assumption (\ref{optimal13.1}), there exists a sufficiently large constant $R_1> 0$ such that it holds
\begin{equation}
\int_{\R^n} u_1(x) \psi_R(x)\, dx >0 \label{t13.6}
\end{equation}
for all $R > R_1$. As a result, combining the estimates from (\ref{t13.2}) to (\ref{t13.6}) gives
\begin{align}
0< \int_{\R^n} u_1(x) \psi_R(x)\, dx &\lesssim J_{R,t}^{\frac{1}{q}} \Big(R^{-2\alpha+ \frac{\alpha+n\beta}{q'}}+ R^{-\alpha- 2\delta_1\beta+ \frac{\alpha+n\beta}{q'}}\Big)+ J_R^{\frac{1}{q}}\, R^{-2\beta+ \frac{\alpha+n\beta}{q'}}- I_R \nonumber \\
&\lesssim J_R^{\frac{1}{q}} \Big(R^{-2\alpha+ \frac{\alpha+n\beta}{q'}}+ R^{-\alpha- 2\delta_1\beta+ \frac{\alpha+n\beta}{q'}}+ R^{-2\beta+ \frac{\alpha+n\beta}{q'}}\Big)- I_R \label{t13.7}
\end{align}
for all $R > R_1$. In the same arguments we may conclude the following estimate for all $R > R_1$:
\begin{align}
0< \int_{\R^n} v_1(x) \psi_R(x)\, dx &\lesssim I_{R,t}^{\frac{1}{p}} \Big(R^{-2\alpha+ \frac{\alpha+n\beta}{p'}}+ R^{-\alpha- 2\delta_2\beta+ \frac{\alpha+n\beta}{p'}}\Big)+ I_R^{\frac{1}{p}}\, R^{-2\beta+ \frac{\alpha+n\beta}{p'}}- J_R \nonumber \\
&\lesssim I_R^{\frac{1}{p}} \Big(R^{-2\alpha+ \frac{\alpha+n\beta}{p'}}+ R^{-\alpha- 2\delta_2\beta+ \frac{\alpha+n\beta}{p'}}+ R^{-2\beta+ \frac{\alpha+n\beta}{p'}}\Big)- J_R. \label{t13.8}
\end{align}
Without loss of generality we can assume $\delta_1 \ge \delta_2$. Now let us fix $\alpha:= 2- 2\delta_1+ \frac{\delta_1- \delta_2}{2(1-\delta_2)}\frac{(nq-n-2q)(n-2)}{1+q}$ and $\beta:= 1- \frac{\delta_1- \delta_2}{2(1-\delta_2)}\frac{nq+2-n}{1+q}$. For this choice, we may verify that
$$ -2\alpha \le -2\beta,\quad -\alpha- 2\delta_1\beta \le -2\beta \quad \text{ and }\quad -\alpha- 2\delta_2\beta \le -2\beta. $$
From (\ref{t13.7}) and (\ref{t13.8}) it follows immediately that
\begin{align*}
I_R&\lesssim J_R^{\frac{1}{q}}R^{-2\beta+ \frac{\alpha+n\beta}{q'}}, \\
J_R&\lesssim I_R^{\frac{1}{p}}R^{-2\beta+ \frac{\alpha+n\beta}{p'}}.
\end{align*}
Therefore, we arrive at
\begin{align}
I_R^{\frac{pq-1}{pq}} &\lesssim R^{- 2\beta+ \frac{\alpha+n\beta}{q'}+ (-2\beta+ \frac{\alpha+n\beta}{p'})\frac{1}{q}}=: R^{\gamma_1}, \label{t13.9} \\ 
J_R^{\frac{pq-1}{pq}} &\lesssim R^{- 2\beta+ \frac{\alpha+n\beta}{p'}+ (-2\beta+ \frac{\alpha+n\beta}{q'})\frac{1}{p}}=: R^{\gamma_2} \label{t13.10}.
\end{align}
It is obvious that the assumption (\ref{optimal13.2}) is equivalent to $\gamma_2 \le 0$. For this reason, we shall divide our attention into two subcases. \medskip

\noindent \textit{Case 1:} Let us consider the subcritical case of $\gamma_2 <0$. Then, we let $R \to \ity$ in (\ref{t13.10}) to obtain
$$ J_R=  \int_0^{\ity}\int_{\R^n}|u(t,x)|^q \eta_R(t,x)\,dxdt= 0, $$
which follows $u \equiv 0$, a contradiction to the assumption (\ref{optimal13.1}). This means that there is no global (in time) Sobolev solution to (\ref{pt1.1}) in the subcritical case. \medskip

\noindent \textit{Case 2:} Let us now come back to the critical case of $\gamma_2= 0$. At first, we introduce the following constants:
\begin{align*}
&C_{u_1}:= \int_{\R^n} u_1(x) \psi_R(x) \quad \text{ and }\quad C_{v_1}:= \int_{\R^n} v_1(x) \psi_R(x), \\ 
&D_{p'}:= \Big(\int_{\R^n} \big< \tilde{x}\big>^{-n-2\delta_0}\, d\tilde{x}\Big)^{\frac{1}{p'}} \quad \text{ and }\quad D_{q'}:= \Big(\int_{\R^n} \big< \tilde{x}\big>^{-n-2\delta_0}\, d\tilde{x}\Big)^{\frac{1}{q'}}.
\end{align*}
After repeating some arguments as we have proved in the subcritical case, we may conclude the following estimates:
\begin{align*}
&0< I_R+ C_{u_1} \le D_{q'} J_R^{\frac{1}{q}}R^{-2\beta+ \frac{\alpha+n\beta}{q'}}, \\ 
&0< J_R+ C_{v_1} \le D_{p'} I_R^{\frac{1}{p}}R^{-2\beta+ \frac{\alpha+n\beta}{p'}}.
\end{align*}
Thus, it follows that
\begin{equation}
J_R+ C_{v_1} \le D_{p'}D^{\frac{1}{p}}_{q'} J_R^{\frac{1}{pq}}R^{\gamma_2}= D_{p'}D^{\frac{1}{p}}_{q'} J_R^{\frac{1}{pq}}. \label{t13.2.11}
\end{equation}
For this reason, we obtain immediately
$$J_R \le D_{p'}D^{\frac{1}{p}}_{q'} J_R^{\frac{1}{pq}} \quad \text{ and }\quad C_{v_1} \le D_{p'}D^{\frac{1}{p}}_{q'} J_R^{\frac{1}{pq}}. $$
Consequently, it implies
\begin{equation}
J_R \le D_0, \label{t13.2.12}
\end{equation}
where $D_0:= \Big(D_{p'}D^{\frac{1}{p}}_{q'}\Big)^{\frac{pq}{pq-1}}$ is a positive constant, and
\begin{equation}
J_R \ge \left(\frac{C_{v_1}}{D_{p'}D^{\frac{1}{p}}_{q'}}\right)^{pq}. \label{t13.2.13}
\end{equation}
By replacing (\ref{t13.2.13}) into the left-hand side of (\ref{t13.2.11}), a direct calculation leads to
$$ J_R \ge \frac{(C_{v_1})^{(pq)^2}}{\Big(D_{p'}D^{\frac{1}{p}}_{q'}\Big)^{pq+(pq)^2}}. $$
Then, we use iteration arguments to arrive at the following estimate for any integer $j\ge 1$:
\begin{equation}
J_R \ge \frac{C_{v_1}^{(pq)^j}}{\Big(D_{p'}D^{\frac{1}{p}}_{q'}\Big)^{pq+(pq)^2+\cdots+(pq)^j}}= \frac{C_{v_1}^{(pq)^j}}{\Big(D_{p'}D^{\frac{1}{p}}_{q'}\Big)^{\frac{(pq)^{j+1}-pq}{pq-1}}}= \Big(D_{p'}D^{\frac{1}{p}}_{q'}\Big)^{\frac{pq}{pq-1}}\left(\frac{C_{v_1}}{\Big(D_{p'}D^{\frac{1}{p}}_{q'}\Big)^{\frac{pq}{pq-1}}}\right)^{(pq)^j}. \label{t13.2.14}
\end{equation}
Let us now choose the constant
$$ \epsilon_2= \int_{\R^n} \big< \tilde{x}\big>^{-n-2\delta_0}\, d\tilde{x} $$
in the assumption (\ref{optimal13.1}). This means that there exists a sufficiently large constant $R_2> 0$ such that
$$ \int_{\R^n} v_1(x) \psi_R(x)\, dx > \epsilon_2 $$
for all $R > R_2$. We can see that the above assumption is equivalent to
$$ C_{v_1}> \int_{\R^n} \big< \tilde{x}\big>^{-n-2\delta_0}\, d\tilde{x}= \Big(D_{p'}D^{\frac{1}{p}}_{q'}\Big)^{\frac{pq}{pq-1}}, \quad \text{that is, }\quad \frac{C_{v_1}}{\Big(D_{p'}D^{\frac{1}{p}}_{q'}\Big)^{\frac{pq}{pq-1}}}> 1. $$
Hence, passing $j \to \ity$ in (\ref{t13.2.14}) gives $J_R \to \ity$. This is a contradiction to the boundedness of $J_R$ in (\ref{t13.2.12}). As a consequence, we may conclude the nonexistence of global (in time) Sobolev solution to (\ref{pt1.1}) in the critical case. Summarizing, the proof of Theorem \ref{dl1.3} is completed.

\subsection{Proof of Theorem \ref{dl1.4}}
We follow the ideas from the proof of Theorem \ref{dl1.3}. We introduce the test functions $\varphi= \varphi(t)$ as in Theorem \ref{dl1.3} and $\psi= \psi(|x|):= \big< x\big>^{-n-2\delta}$. Then, we may repeat exactly, on the one hand, the proof of Theorem \ref{dl1.3} to conclude the following estimates:
\begin{align*}
\int_{\R^n} u_1(x) \psi_R(x)\, dx+ I_R &\le C_{q'} J_R^{\frac{1}{q}} \Big(R^{-2\alpha+ \frac{\alpha+n\beta}{q'}}+ R^{-\alpha- 2\delta\beta+ \frac{\alpha+n\beta}{q'}}+ R^{-2\beta+ \frac{\alpha+n\beta}{q'}}\Big), \\
\int_{\R^n} v_1(x) \psi_R(x)\, dx+ J_R &\le C_{p'} I_R^{\frac{1}{p}} \Big(R^{-2\alpha+ \frac{\alpha+n\beta}{p'}}+ R^{-\alpha- 2\delta\beta+ \frac{\alpha+n\beta}{p'}}+ R^{-2\beta+ \frac{\alpha+n\beta}{p'}}\Big),
\end{align*}
where
$$C_{p'}:= \Big(\int_{\R^n} \big< \tilde{x}\big>^{-n-2\delta}\, d\tilde{x}\Big)^{\frac{1}{p'}} \quad \text{ and }\quad C_{q'}:= \Big(\int_{\R^n} \big< \tilde{x}\big>^{-n-2\delta}\, d\tilde{x}\Big)^{\frac{1}{q'}}. $$
Let us now fix $\alpha:= 2-2\delta$ and $\beta:= 1$. As a result, from the both above estimates we obtain
\begin{align}
\int_{\R^n} u_1(x) \psi_R(x)\, dx+ I_R &\le C_{q'} J_R^{\frac{1}{q}} R^{-2+ \frac{2-2\delta+n}{q'}}, \label{t14.1} \\
\int_{\R^n} v_1(x) \psi_R(x)\, dx+ J_R &\le C_{p'} I_R^{\frac{1}{p}} R^{-2+ \frac{2-2\delta+n}{p'}}, \label{t14.2}
\end{align}
On the other hand, because of the assumption (\ref{optimal14.1}), the following estimate holds:
\begin{align}
\int_{\R^n} u_1(x) \varphi_R(x)\, dx&\ge \epsilon_0 \int_{\R^n}(1+|x|)^{-\frac{n+\epsilon_1}{m}}\, \varphi_R(x)\, dx \nonumber \\ 
&\ge \epsilon_0\,R^n \int_{\R^n} (1+R|\tilde{x}|)^{-\frac{n+\epsilon_1}{m}}\, \varphi(\tilde{x})\, d\tilde{x} \quad \big(\text{by change of variables } \tilde{x}:= R^{-1}x \big) \nonumber \\
&= \epsilon_0\,R^{n-\frac{n+\epsilon_1}{m}} \int_{\R^n} (R^{-1}+|\tilde{x}|)^{-\frac{n+\epsilon_1}{m}}\, (1+|\tilde{x}|^2)^{-\frac{n+2\delta}{2}}\, d\tilde{x} \ge C_1\epsilon_0 R^{n-\frac{n+\epsilon_1}{m}} \label{t14.3}
\end{align}
for all $R > R_0$, where $R_0> 0$ is a sufficiently large number and $C_1$ is a suitable positive constant. In the same way we also derive
\begin{equation}
\int_{\R^n} v_1(x) \varphi_R(x)\, dx\ge C_2\epsilon_0 R^{n-\frac{n+\epsilon_2}{m}} \label{t14.4}
\end{equation}
for all $R > R_0$, where $C_2$ is a suitable positive constant. Combining the estimates from (\ref{t14.1}) to (\ref{t14.4}) we may arrive at
\begin{align}
C_1\epsilon_0 R^{n-\frac{n+\epsilon_1}{m}}&\le C_{q'}C_{p'}^{\frac{1}{q}} I_R^{\frac{1}{pq}}\,R^{- 2+ \frac{2-2\delta+n}{q'}+ (-2+ \frac{2-2\delta+n}{p'})\frac{1}{q}}- I_R, \label{t14.5} \\
C_2\epsilon_0 R^{n-\frac{n+\epsilon_2}{m}}&\le C_{p'}C_{q'}^{\frac{1}{p}} J_R^{\frac{1}{pq}}\,R^{- 2+ \frac{2-2\delta+n}{p'}+ (-2+ \frac{2-2\delta+n}{q'})\frac{1}{p}}- J_R, \label{t14.6}
\end{align}
for all $R > R_0$. Moreover, applying the inequality
$$ A\,y^\gamma- y \le A^{\frac{1}{1-\gamma}} \quad \text{ for any } A>0,\, y \ge 0 \text{ and } 0< \gamma< 1 $$
to (\ref{t14.6}) leads to
$$ C_2\epsilon_0 R^{n-\frac{n+\epsilon_2}{m}}\le \Big(C_{p'}C_{q'}^{\frac{1}{p}}R^{- 2+ \frac{2-2\delta+n}{p'}+ (-2+ \frac{2-2\delta+n}{q'})\frac{1}{p}}\Big)^{\frac{pq}{pq-1}}= \Big(C_{p'}C_{q'}^{\frac{1}{p}}\Big)^{\frac{pq}{pq-1}}R^{n- 2\delta- \frac{2(1+q)}{pq}} $$
for all $R > R_0$. It follows immediately
\begin{equation}
\frac{C_2\epsilon_0}{\Big(C_{p'}C_{q'}^{\frac{1}{p}}\Big)^{\frac{pq}{pq-1}}}\le R^{- 2\delta- \frac{2(1+q)}{pq}+\frac{n+\epsilon_2}{m}} \label{t14.7}
\end{equation}
for all $R > R_0$. Without loss of generality we can assume $p\le q$, it is clear that the assumption (\ref{optimal14.2}) is equivalent to
$$ \frac{n-2m\delta}{2m}< \frac{1+ q}{pq-1}, $$
that is, $- 2\delta- \frac{2(1+q)}{pq}+\frac{n}{m}< 0$. Then, we can choose a sufficiently small constant $\e>0$ such that the following relation still holds $- 2\delta- \frac{2(1+q)}{pq}+\frac{n+\e}{m}< 0$. Now we take $\epsilon_2=\e$ in the assumption (\ref{optimal14.1}). By letting $R \to \ity$ in (\ref{t14.7}) we obtain a contradiction to the choice of positive constants $\epsilon_0,\,C_2,\,C_{p'}$ and $C_{q'}$. Summarizing, the proof of Theorem \ref{dl1.4} is completed.

\section*{Acknowledgment}
The PhD study of MSc. T.A. Dao is supported by Vietnamese Government's Scholarship (Grant number: 2015/911). The author would like to thank sincerely to Prof. Michael Reissig for valuable discussions and Institute of Applied Analysis for their hospitality. The author is grateful to the referee for his careful reading of the manuscript and for helpful comments.

\section*{Appendix}

\begin{md}[Fractional Gagliardo-Nirenberg inequality] \label{fractionalGagliardoNirenberg}
Let $1<p,\, p_0,\, p_1<\infty$, $\sigma >0$ and $s\in [0,\sigma)$. Then, it holds the following fractional Gagliardo-Nirenberg inequality for all $u\in L^{p_0} \cap \dot{H}^\sigma_{p_1}$:
$$ \|u\|_{\dot{H}^{s}_p}\lesssim \|u\|_{L^{p_0}}^{1-\theta}\,\, \|u\|_{\dot{H}^{\sigma}_{p_1}}^\theta, $$
where $\theta=\theta_{s,\sigma}(p,p_0,p_1)=\frac{\frac{1}{p_0}-\frac{1}{p}+\frac{s}{n}}{\frac{1}{p_0}-\frac{1}{p_1}+\frac{\sigma}{n}}$ and $\frac{s}{\sigma}\leq \theta\leq 1$.
\end{md}
For the proof one can see \cite{Ozawa}.


\end{document}